\documentclass[a4paper,12pt]{amsart}

\newtheorem{theo}{Theorem}[section]

\newtheorem{prop}[theo]{Proposition}
\newtheorem{lem}[theo]{Lemma}
\newtheorem{cor}[theo]{Corollary}

\newtheorem{rema}[theo]{Remark}
\newtheorem{remas}[theo]{Remarks}

\def \kbar {{\overline k}}
\def \Xbar {{\overline X}}

\def \Romannumeral #1 {\expandafter\uppercase\expandafter {\romannumeral #1} }
\def \br {{\rm{Br}}}
\def \brnr {{\rm{Br_{\rm nr}}}}

\def \P {{\bf P}}

\def \pic {{\rm {Pic}}}
\def \div {{\rm{Div}}}
\def \gal {{\rm{Gal}}}
\def \Alb {{\rm{Alb}}}

\def \calo {{\mathcal O}}
\def \spec {{\rm{Spec\,}}}

\def \Hom {{\rm {Hom}}}
\def \Ext {{\rm {Ext}}}

\def \Z {{\bf Z}}
\def \Q {{\bf Q}}

\def \Alb {{\rm Alb}}

\def \cok {{\rm{coker\,}}}
\def \im {{\rm {Im\,}}}
\def \G {{\bf G}_m}
\def \A {{\bf A}}
\def \Het {H_{\mbox{\scriptsize\'et}}}

\def\smallsquare{\vbox{\hrule\hbox{\vrule height 1 ex\kern 1 ex\vrule}\hrule}}
\def\enddem{\hfill \smallsquare\vskip 3mm}

\usepackage{amscd}
\usepackage{amssymb}
\usepackage{amsmath}
\usepackage{pb-diagram}

\everymath{\displaystyle}

\def \hyp {{\bf H}}

\DeclareFontFamily{U}{wncy}{}
\DeclareFontShape{U}{wncy}{m}{n}{%
   <5>wncyr5%
   <6>wncyr6%
   <7>wncyr7%
   <8>wncyr8%
   <9>wncyr9%
   <10>wncyr10%
   <11>wncyr10%
   <12>wncyr6%
   <14>wncyr7%
   <17>wncyr8%
   <20>wncyr10%
   <25>wncyr10}{}
\DeclareMathAlphabet{\cyrille}{U}{wncy}{m}{n}
\def\Sha{\cyrille X}
\def\Be{\cyrille B}
\def\be{\cyrille b}

\title{Local-global principles for 1-motives}
\author{David Harari and Tam\'as Szamuely}
\address{Universit\'e de Paris-Sud Math\'ematique, B\^atiment 425, 91405 Orsay, France}
\email{David.Harari@math.u-psud.fr}
\address{Alfr\'ed R\'enyi Institute of Mathematics, Hungarian Academy of Sciences, PO Box 127, H-1364 Budapest, Hungary}
\email{szamuely@renyi.hu}

\begin{document}
\begin{abstract} Building upon our arithmetic duality theorems for 1-motives, we prove
that the Manin obstruction related to a finite subquotient $\Be (X)$ of the Brauer group
is the only obstruction to the Hasse principle for rational points on torsors under
semiabelian varieties over a number field, assuming the finiteness of the Tate-Shaferevich
group of the abelian quotient. This theorem answers a question by Skorobogatov in the
semiabelian case, and is a key ingredient of recent work on the elementary obstruction for
homogeneous spaces over number fields. We also establish a Cassels--Tate type dual exact
sequence for 1-motives, and give an application to weak approximation.\end{abstract}
\maketitle

\section{Introduction}

In this paper we use the duality theorems of \cite{hasza} to prove some results related to
1-motives over number fields. The main object of study is the Manin obstruction to the
Hasse principle and weak approximation on torsors under a semi-abelian variety over a
number field $k$ (i.e. a commutative $k$-group scheme which is an extension of an abelian
variety by a torus).

We briefly recall the basic idea of the Manin obstruction. Given a smooth, geometrically
integral variety $X$ over our number field $k$, one considers the ring of adeles $\A_k$ of
$k$, and defines a pairing
\begin{equation}\label{pairing}
X({\bf A}_k)\times \br\, X\to {\bf Q}/{\bf Z}
\end{equation}
by evaluating elements of the cohomological Brauer group $\br\, X:=\Het^2(X, \G)$ at each
component and then taking the sum of local invariants (which is known to be finite; see
e.g. \cite{skobook}, p. 101). By global class field theory, the diagonal image of $X(k)$
in $X({\bf A}_k)$ is  contained in the subset $X({\bf A}_k)^{\br}$ of adeles annihilated
by the above pairing. Consequently, the emptiness of $X({\bf A}_k)^{\br}$ is an
obstruction to the Hasse principle if $X({\bf A}_k)$ itself is nonempty.

In our case, a special role is played by a subquotient $\Be(X)$ of $\br X $ defined as
follows. Fix an algebraic closure $\kbar$ of $k$. First one defines the subgroup $\br_1\,
X\subset \br X $ as the kernel of the natural map $\br X\to \br (X\times_k\kbar)$, and
then defines $\br_a\,X$ as the quotient of $\br_1\,X$ modulo the image of the map $\br\,
k\to \br\, X$. Finally one puts
$$
\Be (X):= {\rm Ker} (\br_a\,X\to \prod_{v\in\Omega}\br_a (X\times_kk_v)),
$$
where $\Omega$ denotes the set of all places of $k$. (Note that some authors use the
notation $\Be (X)$ for the preimage of our $\Be (X)$ in $\br\, X$.) The pairing
(\ref{pairing}) manifestly factors through the image of $\br\, k$ in $\br X$, and hence
induces a pairing
\begin{equation}\label{bepairing}
X({\bf A}_k)\times \Be (X)\to {\bf Q}/{\bf Z}.
\end{equation}

Our first main result concerns a semi-abelian variety $G$ over $k$, i.e. an extension of
an abelian variety $A$ by a torus $T$. Recall that the {\it Tate-Shafarevich} group of $A$
is defined as the kernel

$$\Sha(A):=\ker\, (H^1(k,A) \to \prod_{v \in \Omega} H^1(k_v,A))$$
of the natural restriction map in Galois cohomology.

\begin{theo}\label{thm1}
Let $G$ be a semi-abelian variety defined over $k$, and let $X$ be a $k$-torsor under $G$.
Assume that the Tate-Shafarevich group of the abelian quotient $A$ of $G$ is finite. If
there is an adelic point of $X$ annihilated by all elements of $\Be (X)$ under the pairing
(\ref{bepairing}), then $X$ has a $k$-rational point.
\end{theo}
The theorem answers positively a question raised by Skorobogatov in \cite{skobook} (p.
133, question 1) for semi-abelian varieties. The analogous result for connected linear
algebraic groups has been known for a long time (see \cite{sansuc}). However, recently
Borovoi, Colliot-Th\'el\`ene and Skorobogatov (\cite{bcs}, Proposition 3.16 together with
Theorem 2.12) gave an example of a connected non-commutative and non-linear algebraic
group over $\Q$ for which the statement fails.

In the same paper Borovoi, Colliot-Th\'el\`ene and Skorobogatov  also prove a positive
result that relies on our Theorem \ref{thm1}. Namely, they show that if $k$ is a totally
imaginary number field, the existence of a $k$-point on a torsor $X$ under a connected
algebraic group satisfying a finiteness condition on the abelian quotient as in Theorem
\ref{thm1} is equivalent to the existence of a Galois-equivariant splitting of the exact
sequence
\begin{equation}\label{elob}
0\to \kbar^\times\to \kbar(X)^\times\to \kbar(X)^\times/\kbar^\times \to 0
\end{equation}
 where
$\kbar(X)^\times$ is the group of invertible rational functions on $X\times_k\kbar$. To
obtain the statement they really have to work with $\Be(X)$, the weaker version with the
pairing (\ref{pairing}) is not sufficient. Another nice application having a similar
flavour is the following recent result of O. Wittenberg \cite{witten}: assuming the
finiteness of the Tate-Shafarevich groups of abelian varieties over $k$, the sequence
(\ref{elob}) splits for an arbitrary smooth $k$-variety $X$ if and only if the pairing
(\ref{bepairing}) is trivial. Wittenberg uses our Theorem \ref{thm1} via a result in
\cite{esch}.

Another reason why it is more interesting to work with the subquotient $\Be (X)$ instead
of the whole group $\br \, X$ is that under the assumptions of the theorem it is finite
(see Remark \ref{berema} below). In this respect Theorem \ref{thm1} improves the main
result of the {\em refuznik} paper \cite{imrn} that shows that a statement as in Theorem
\ref{thm1} holds even for arbitrary connected algebraic groups, provided that one replaces
$\Be (X)$ with  the unramified part of $\br_1 X$, which is a much bigger group.

Though Theorem \ref{thm1} has been known for quite some time in the extreme cases $G=A$
(Manin \cite{manin}) and $G=T$ (Sansuc \cite{sansuc}), the general case does not follow
from these, and is substantially more difficult. Our approach is based on a strategy
similar to the proofs in the extreme cases, but it is more conceptual and avoids some
rather involved cocycle calculations that made earlier texts hard to follow (at least for
us). Among new ingredients we use the theory of generalised Albanese varieties and
1-motives, and a new interpretation of the pairing (\ref{bepairing}) as a cup-product in
\'etale hypercohomology. The latter fact is valid for an arbitrary smooth variety and is
of independent interest (see Section \ref{secbm}).

\medskip

To state  our second main result, let $M=[Y\to G]$ be a 1-motive over $k$ (i.e. a complex
of $k$-group schemes placed in degrees $(-1,0)$ with $Y$ \'etale locally isomorphic to
$\Z^r$ for some $r\geq 0$ and $G$ a semiabelian variety). Denote the dual 1-motive
(\cite{deligne}, 10.2.1) by $M^*=[Y^*\to G^*]$. For each positive integer $i$ denote by
$\Sha^i(M)$ (resp. $\Sha^i_\omega(M)$) the subgroup of $\hyp^i(k,M)$ consisting of those
elements whose restriction to $\hyp^i(k_v,M)$ is zero for all places (resp. for all but
finitely many places) of $k$. In Section \ref{catate} we shall prove the following
generalisation of the classical Cassels-Tate dual exact sequence to 1-motives.

\begin{theo} \label{casstatetheo} Assume that the Tate-Shafarevich group $\Sha
  (A)$  of the abelian quotient of $G$ is finite.
Then there is an exact sequence of abelian groups
$$0 \to \overline{{\hyp^0(k,M)}}\to \prod_{v\in\Omega} \hyp^0(k_v,M) \to \Sha^1_{\omega}(M^*)^D
\to \Sha^1(M) \to 0,$$ where $\overline{{\hyp^0(k,M)}}$ denotes the closure of the
diagonal image of ${{\hyp^0(k,M)}}$ in the topological product  of the $\hyp^0(k_v,M)$,
and $A^D:=\Hom(A,\Q/\Z)$ for a discrete abelian group $A$. By convention, for $v$
archimedean we take here the modified (Tate) hypercohomology groups instead of the usual
ones.
\end{theo}

The third (resp. fourth) maps in the above sequence are induced by the local (resp.
global) duality pairings of \cite{hasza} (see Section \ref{catate} for more details), and
the topology on $\hyp^0(k_v,M)$ was defined in \S 2 of the same reference. A similar
statement for connected linear algebraic groups was proven by Sansuc (\cite{sansuc},
Theorem 8.12). Our method yields a new proof of his result in the case of tori.

The case $M=[0\to G]$ of this exact sequence allows us to give a new short proof of the
weak approximation part of the main result of \cite{imrn} in the crucial case of a
semi-abelian variety. For the precise statement, see Section \ref{secaf}.

We thank Klaus K\"unnemann and Olivier Wittenberg for helpful discussions, and the referee
for his or her very careful reading of the manuscript. The second author acknowledges
support from OTKA grant No. 61116 as well as from the BUDALGGEO Intra-European project,
and is grateful to the Institut Henri Poincar\'e for its hospitality.

\section{Preliminaries on the Brauer group}

As a preparation for the proof of Theorem \ref{thm1}, we collect here some auxiliary
statements concerning the Brauer group. Statements \ref{lemupic}--\ref{braremas} will not
serve until Section \ref{secthm1}.

We investigate the exact sequence of complexes of $\gal(\kbar|k)$-modules
\begin{equation}\label{extA}
0\to [\kbar^\times\to 0]\to [\kbar (X)^\times\to \div \, \Xbar]\to
[\kbar(X)^\times/\kbar^\times\to \div\, \Xbar]\to 0
\end{equation}
for an arbitrary smooth geometrically integral variety $X$ over a field $k$. Here, as
usual, $\Xbar$ stands for $X\times_k\kbar$, $\div\,\Xbar$ for the group of divisors on
$\Xbar$ and $\kbar(X)$ for its function field. In accordance with our conventions for
1-motives, the above complexes are placed in degrees $-1$ and 0.

\begin{lem}\label{bra}
There are canonical isomorphisms
$$
\hyp^1(k, [\kbar (X)^\times\to \div \, \Xbar])\cong \br_1\, X$$ and, assuming $H^3(k,
\kbar^\times)=0$, $$ \hyp^1(k, [\kbar (X)^\times/\kbar^\times\to \div \, \Xbar])\cong
\br_a\, X.
$$
\end{lem}

\begin{dem} For the first isomorphism we use the long exact hypercohomology sequence associated with the distinguished triangle
$$
\kbar (X)^\times\to \div \, \Xbar\to [\kbar (X)^\times\to \div \, \Xbar]\to
\kbar(X)^\times [1]
$$
and the fact that the permutation module $\div\,\Xbar$ has trivial $H^1$ to obtain
$$
\hyp^1(k, [\kbar (X)^\times\to \div \, \Xbar])\cong\ker\left(H^2(k, \kbar(X)^\times)\to
H^2(k, \div\,\Xbar)\right).
$$
The latter group is identified in \cite{skobook} with $\br_1 X$: see the second column of
the exact diagram (4.17) on p. 72. The injectivity of the map $\br_1\, X\to
\ker\left(H^2(k, \kbar(X)^\times)\to H^2(k, \div\,\Xbar)\right)$ constructed there follows
again from the vanishing of $H^1(k, \div\, \Xbar)$, and the assumption
$\kbar[X]^\times=\kbar^\times$ made in the reference is not used at this point.

Another (?) argument is to use the quasi-isomorphism of Galois modules $[\kbar
(X)^\times\to \div \, \Xbar]\stackrel\sim\to \tau_{\leq 1}{\bf R}\pi_*\G[1]$ constructed
in \cite{bvh}, Lemma 2.3, where $\pi:\, \Xbar\to \spec\kbar$ is the natural projection. It
yields an isomorphism of $\hyp^1(k, [\kbar (X)^\times\to \div \, \Xbar])$ with $\hyp^2(k,
\tau_{\leq 1}{\bf R}\pi_*\G)$, in turn isomorphic to $\ker(H^2(X, \G)\to H^0(k, H^2(\Xbar,
\G))$ by a Hochschild-Serre argument.

The second isomorphism of the lemma follows from the first one in view of the long exact
cohomology sequence associated with (\ref{extA}) and the assumption $H^3(k,
\kbar^\times)=0$.\end{dem}

The complex of Galois modules $[\kbar (X)^\times/\kbar^\times\to \div \, \Xbar]$ was
considered in several recent papers, in particular in Borovoi and van Hamel \cite{bvh} and
Colliot-Th\'el\`ene \cite{flasque}. The following property seems to have been observed by
all of us:

\begin{lem}\label{lemupic} Assume that $X$ has a smooth compactification $X^c$ over $k$. Then there is a natural Galois-equivariant quasi-isomorphism of complexes
$$[\kbar (X)^\times/\kbar^\times\to \div \, \Xbar]\simeq [\div^\infty\, \Xbar^c\to\pic\,\Xbar^c],
$$
where $\div^\infty\, \Xbar^c$ denotes the group of divisors on $\Xbar^c:=X^c\times_k\kbar$
supported in $\Xbar^c\setminus \Xbar$, and the maps in both complexes are divisor maps.
\end{lem}

\begin{dem} We consider the map of two-term complexes
\begin{equation}\label{upicdiag}
\begin{CD}
[\kbar (X)^\times/\kbar^\times @>>> \div \, \Xbar] \\
@VVV @VVV \\
[\div^\infty\, \Xbar^c @>>> \pic\,\Xbar^c]
\end{CD}
\end{equation}
where the right vertical map is induced by sending a codimension 1 point on $\Xbar$ to the
class of the corresponding point of $\Xbar^c$ {\em with a minus sign.} This sign
convention implies that we indeed have a map of complexes, and it is a quasi-isomorphism
by construction.
\end{dem}

A combination of the two previous lemmas yields:

\begin{cor}\label{corbra} If moreover $k$ satisfies $H^3(k, \kbar^\times)=0$, we have a natural isomorphism
$$
\br_a \, X\cong\hyp^1(k, [\div^\infty\, \Xbar^c\to\pic\,\Xbar^c]).
$$
\end{cor}

If $k$ is of characteristic 0, the smooth compactification $X^c$ always exists thanks to
Hironaka's theorem on resolution of singularities. In particular, the statement of the
corollary holds over number fields and their completions, since they are of characteristic
0 and satisfy the cohomological condition.

\begin{remas}\label{braremas}\rm ${}$
\begin{enumerate}
\item Note that when $X$ is proper, the isomorphism of the corollary is just  the well-known identification $\br_a\, X\cong H^1(k, \pic \,\Xbar)$ induced by the Hochschild-Serre spectral sequence.
\item We shall
also need a sheafified version of Lemma \ref{lemupic} over sufficiently small nonvoid open
subsets of $U\subset\spec\calo_k$, where $\calo_k$ denotes the ring of integers of $k$, as
usual. For $U$ sufficiently small the $k$-varieties $X$ and $X^c$ extend to smooth
$U$-schemes ${\mathcal X}$ and ${\mathcal X}^c$, with ${\mathcal X}^c$ projective over
$U$, and ${\mathcal X}$ open in ${\mathcal X}^c$. We shall work on the {\em big \'etale
site} of $U$ restricted to the subcategory $Sm/U$ of smooth $U$-schemes of finite type.
Consider the \'etale sheaf $\div_{{\mathcal X}^c/U}$ associated with the presheaf
$S\mapsto \div({\mathcal X}^c\times_U S/S)$ of relative Cartier divisors on this site (see
e.g. \cite{blr} pp. 212--213, for a discussion of relative effective Cartier divisors, and
then take group completion). It contains as subsheaves the sheaf $\div_{{\mathcal X}/U}$
of relative divisors on $\mathcal X\times_U S/S$ and the sheaf $\div^{\infty}_{{\mathcal
X}^c/U}$ of relative divisors on ${\mathcal X}^c\times_US/S$ with support in $({\mathcal
X}^c \setminus {\mathcal X})\times_US$. There is also the sheaf $\pic_{{\mathcal X}^c/U}$
given by the relative Picard functor. Finally, denote by ${\mathcal K}_{\mathcal
X}^\times$ the \'etale sheaf on $Sm/U$ associated with the presheaf sending $S$ to the
group of rational functions on $\mathcal X\times_U S$ that are regular and invertible on a
dense open subset of each fibre of the projection onto $S$. It contains as a subsheaf the
pullbacks of invertible functions on $S$; we identify it with the sheaf $\G$ on the big
\'etale site of $Sm/U$. We contend that we may lift the quasi-isomorphism (\ref{upicdiag})
to a quasi-isomorphism
$$ \begin{CD}
[{\mathcal K}_{\mathcal X}^\times/\G @>>> \div_{{\mathcal X}/U}]  \cr @VVV @VVV \cr
[\div^{\infty}_ {{\mathcal X}^c/U} @>>> \pic_{{\mathcal X}^c/U}]
\end{CD} $$ of
complexes of \'etale sheaves on $Sm/U$. The definition of the morphism is as above (one
works with the inclusion ${\mathcal X}\times_US\subset {\mathcal X}^c\times_US$ for each
object $S$ of $Sm/U$). To check that it is a quasi-isomorphism, we may restrict to the
small \'etale site of each $S$. Then it is sufficient to check it at every geometric point
$\bar u$ of  $S$, where the claim follows from Lemma \ref{lemupic} (which is valid in any
characteristic).
\end{enumerate}
\end{remas}

\section{Reinterpretations of the Brauer-Manin pairing}\label{secbm}

In this section we use exact sequence (\ref{extA}) and Lemma \ref{bra} to give other
formulations of the Brauer-Manin pairing
$$
X(\A_k)\times \Be(X)\to \Q/\Z.
$$
Our $X$ is still an arbitrary smooth geometrically integral variety over a number field
$k$, and we assume $X(\A_k)\neq \emptyset$.

We start with a couple of well-known observations. Since the elements of $\Be(X)$ are
locally constant by definition, the maps
\begin{equation}\label{bemap}
\be:\,\Be(X)\to\Q/\Z
\end{equation}
given by evaluation on an adelic point $(P_v)$ do not depend on the choice of $(P_v)$, so
defining the pairing is equivalent to defining the map $\be$. There is a commutative
diagram with exact rows
\begin{equation}\label{braseq}
\begin{CD}
0 @>>> \br \, k @>>> \br_1\, X @>>> \br_a\, X @>>> 0 \\
&& @VVV @VVV @VVV \\
0 @>>>  \bigoplus_{v\in\Omega} \br\, k_v @>>> \bigoplus_{v\in\Omega} \br_1 \,X_v @>>>
\bigoplus_{v\in\Omega} \br_a \,X_v @>>> 0
\end{CD}
\end{equation}
where $X_v:= X\times_kk_v$, and the first map in the bottom row is injective because our
assumption that $X(\A_k)\neq \emptyset$ implies the injectivity of each map $\br\, k_v\to
\br_1(X\times_k{k_v})$. The first map in the top row is then injective because so is the
left vertical map, by the Hasse principle for Brauer groups. Applying the snake lemma to
the diagram we thus have a map
$$
\Be(X)=\ker(\br_a\, X \to \bigoplus_{v\in\Omega} \br_a \, X_v) \to\cok(\br\, k \to
\bigoplus_{v\in\Omega} \br\, {k_v})\cong \Q/\Z.
$$

\begin{lem}\label{belem}
The above map equals the map $\be:\, \Be(X)\to\Q/\Z$.
\end{lem}

\begin{dem} For $\alpha\in\Be(X)$ the value $\be(\alpha)$ is defined by lifting first $\alpha$ to $\alpha'\in\br_1\, X$, then sending $\alpha'$ to an element of $\oplus_v \br\, k_v$ via a family of local sections $(s_v:\, \br_1 X\to \br\, k_v)$ determined by an adelic point of $X$, and finally taking the sum of local invariants. Since each $s_v$ factors through $\br_1 (X\times_k{k_v})$, this yields the same element as the snake lemma construction.
\end{dem}

Now observe that in view of Lemma \ref{bra} one may also obtain the diagram (\ref{braseq})
by taking the long exact hypercohomology sequence coming from the diagram

\vbox{\small
$$0\to [\kbar^\times\to 0]\to \qquad [\kbar (X)^\times\to \div \, \Xbar]\quad\to\quad
[\kbar(X)^\times/\kbar^\times\to \div\, \Xbar]\to 0
$$
$$
\downarrow\qquad\quad\qquad\qquad \qquad\qquad \downarrow\qquad\qquad\qquad\qquad
\qquad\downarrow
$$
$$0\to \bigoplus_{v\in\Omega}[\kbar_v^\times\to 0]\to\bigoplus_{v\in\Omega} [\kbar_v (X)^\times\to \div \, \Xbar_v]\to \bigoplus_{v\in\Omega}[\kbar_v(X_v)^\times/\kbar_v^\times\to \div\, \Xbar_v]\to 0,
$$}

\noindent where $\Xbar_v:=X\times_k\kbar_v$. The zeros on the right in (\ref{braseq}) come
from the fact that the groups $H^3(k, \G)$ and $H^3(k_v, \G)$ all vanish.

\begin{rema}\label{secrema}\rm Note in passing that the sections $s_v$ used in the above proof come from  Galois-equivariant splittings
\begin{equation}\label{split}
[\kbar_v (X)^\times\to \div (X\times_k\kbar_v)]\to [\kbar_v^\times\to 0]
\end{equation}
of the base change of the extension (\ref{extA}) to $\kbar_v$. As maps of complexes, the
latter are given in degree $-1$ by a natural splitting of the inclusion map
$\kbar_v^\times\to\kbar_v(X)^\times$ coming from $P_v$ as constructed e.g. in
(\cite{skobook}, Theorem 2.3.4 (b)), and in degree 0 by the zero map. In particular, the
extension (\ref{extA}) is locally split.
\end{rema}

As in Remark \ref{braremas} (2), we now pass to sheaves over the \'etale site of $Sm/U$,
where $U\subset\spec\calo_k$ is a suitable open subset. We can then extend the upper row
of the last diagram to an exact sequence
\begin{equation}\label{KDseq}
0\to \G[1]\to \mathcal{KD(X)}\to \mathcal{KD'(X)}\to 0
\end{equation}
of complexes of \'etale sheaves on $Sm/U$, where
$$
\mathcal{KD(X)}:=[{\mathcal K}_{\mathcal X}^\times\to \div_{{\mathcal X}/U}]$$ and $$
\mathcal{KD'(X)}:=[{\mathcal K}_{\mathcal X}^\times/\G\to \div_{{\mathcal X}/U}].
$$
By Lemma \ref{lemupic} we have $\br_a\, X\cong H^1(k, [\kbar(X)^\times/\kbar^\times\to
\div\, \Xbar])$, so each element of $\br_a\, X$ comes from an element
in $\hyp^1(U, \mathcal{KD'(X)})$. Now assume moreover  $\alpha\in\br_a\, X$ is
locally trivial, i.e. lies in $\Be(X)$.  For a finite place $v$ of $k$ we
have $\hyp^1(k_v^h,
j_v^{h*}\mathcal{KD'(X)})\cong \hyp^1(k_v,
j _v^*\mathcal{KD'(X)})$, where $k_v^h$ is the henselisation of $k$ at $v$,
and ${j_v^h:\spec k_v^h\to U}$ as well as  ${j_v:\spec k_v\to U}$ are the
natural maps. This is shown using the quasi-isomorphism of Lemma \ref{lemupic}, and
then reasoning as in the proof of
(\cite{hasza}, Lemma 2.7).  Next recall that in the case when $k$ is totally
imaginary the arithmetic compact support hypercohomology $\hyp^1_c(U, \mathcal
F)$ of  a complex of sheaves $\mathcal F$ is defined  by $\hyp^i_c(\spec\calo_k,
j_!\mathcal{F})$, where $j:\, U\to\spec\calo_k$ is the natural inclusion. It fits into a
long exact sequence
$$
\dots\to \hyp^1_c(U, \mathcal{F})\to \hyp^1(U, \mathcal{F})\to
\bigoplus_{v\in\spec\calo_k\setminus U}\hyp^1(k_v^h, j_v^{h*}\mathcal{F}) \to\dots
$$
In the general case there are corrective terms coming from the real places;
see the discussion at the beginning of \S 3 in \cite{hasza} (but note the misprint in
formula (8) there: the $\hat k_v$ should be $k_v$ in that paper's notation).
It then follows from the above discussion that we
may lift $\alpha\in\Be(X)$ to an element $\alpha_U\in \hyp^1_c(U, \mathcal{KD'(X)})$ for sufficiently
small $U$. 

There is a cup-product pairing
$$
\cup:\,\Ext^1_{Sm/U}(\mathcal{KD'(X)}, \G[1])\times \hyp^1_c(U, \mathcal{KD'(X)})\to
H^3_c(U, \G)\cong \Q/\Z,
$$
where the Ext-group is taken in the category of \'etale sheaves on $Sm/U$, and the last
isomorphism comes from global class field theory (see \cite{adt}, p. 159). We shall be
interested in the class $\mathcal{E}_X\cup \alpha_U$, where ${\mathcal E}_X$ is the class
of the sequence (\ref{KDseq} )in $\Ext^1_{Sm/U}(\mathcal{KD'(X)}, \G[1])$. Note that there
is a commutative diagram

\vbox{$$ \Ext^1_{Sm/U}(\mathcal{KD'(X)}, \G[1])\times \hyp^1_c(U, \mathcal{KD'(X)})\to
H^3_c(U, \G)\cong \Q/\Z,
$$

$$
\downarrow\qquad\qquad\qquad \qquad\qquad \downarrow\cong\quad\quad\qquad\qquad
\quad\downarrow\, {\rm id}
$$

$$
\Ext^1_{U}(g_*\mathcal{KD'(X)},\G[1])\times \hyp^1_c(U, g_*\mathcal{KD'(X)})\to H^3_c(U,
\G)\cong \Q/\Z,
$$}

\noindent where $g_*$ is the natural pushforward (or restriction) map from the \'etale
site of $Sm/U$ to the small \'etale site of $U$, and the Ext-group in the bottom row is an
Ext-group for \'etale sheaves on $U$. The left vertical map exists because the functor
$g_*$ is exact (and $\G$ as an \'etale sheaf on $U$ is the pushforward by $g$ of the $\G$
on $Sm/U$). The middle isomorphism comes from the fact that the hypercohomology of
complexes of sheaves on the big \'etale site of $U$ equals the hypercohomology on the
small \'etale site. So instead of ${\mathcal E}_X$ we may work with its image
$g_*{\mathcal E}_X$ in $\Ext^1_{U}(g_*\mathcal{KD'(X)},\G[1])$, and omit the $g_*$ from
the notation when no confusion is possible. Note that the generic stalk of $g_*{\mathcal
E}_X$ is the class of the extension (\ref{extA}) in the group
$\Ext^1_k([\kbar(X)^\times/\kbar^\times\to \div\, \Xbar], \kbar^\times[1])$.

\begin{prop}\label{BM} With notations as above, we have $$\mathcal{E}_X\cup \alpha_U=\be(\alpha).$$
\end{prop}

Before starting the proof, note that though one has several choices for $\alpha_U$, the
cup-product depends only on $\alpha$. Indeed two choices of $\alpha_U$ differ by an
element of the direct sum of the groups $\hyp^0(k_v^h, \mathcal{KD'(X)})$, and the
cup-product of each such group with $\Ext^1_U(\mathcal{KD'(X)}, \G[1])$ factors through
the cup-product with $\Ext^1_{k_v^h}(j_v^*\mathcal{KD'(X)}, \G[1])$. But the image of the
class $g_*\mathcal{E}_X$ in these groups is 0, because the extension is locally split
(Remark \ref{secrema}).\medskip

\begin{dem} In order to avoid complicated notation, we do the verification in the case
when $U$ is totally imaginary and the simpler definition of compact support cohomology is
available, and leave the general case to anxious readers. We may work over the small
\'etale site of $U$ by the previous observations; in particular, we identify the complexes
$\mathcal{KD(X)}$ and $\mathcal{KD'(X)}$ with their images under $g_*$.

The cup-product $\mathcal{E}_X\cup \alpha_U$ is none but the image of $\alpha$ by the
boundary map $\hyp^1_c(U, \mathcal{KD'(X)})\to H^3_c(U, \G)$ coming from the long exact
hypercohomology sequence associated with (\ref{KDseq}). Now consider the commutative exact
diagram

\vbox{$$ 0 \quad\to\quad \G[1]\qquad \to\qquad \mathcal{KD(X)}\qquad \to\qquad
\mathcal{KD'(X)} \quad\to\quad 0
$$

$$
\downarrow\quad\qquad\qquad \qquad\qquad \downarrow\qquad\qquad\qquad\qquad\quad\downarrow
$$

$$
0 \to \bigoplus_{v\notin U}j_{v*}j_v^*\G[1]\to \bigoplus_{v\notin
U}j_{v*}j_v^*\mathcal{KD(X)} \to \bigoplus_{v\notin U}j_{v*}j_v^*\mathcal{KD'(X)} \to 0
$$}

\noindent of complexes of \'etale sheaves on $U$, and denote the cones of the vertical
maps by $\mathcal C$, $\mathcal{C_K}$ and $\mathcal{C_K'}$, respectively. The group
$\hyp^1(U, \mathcal{C_K'}[-1])$ may be identified with $\hyp^1_c(U, \mathcal{KD'(X)})$
(apply, for instance, \cite{adt}, Lemma II.2.4 and its proof with our $U$ as $V$ and our
$\spec\calo_k$ as $U$), so we may view $\alpha_U$ as an element of the former group.   The
cup-product $\mathcal{E}_X\cup \alpha_U$ maps to a class in $\hyp^2(U,
\mathcal{C}[-1])=\hyp^1(U, \mathcal{C})$. But when one makes $U$ smaller and smaller and
passes to the limit, this class yields an element in the cokernel of the map $\br\, k \to
\oplus_{v} \br\, {k_v^h}$ which (noting the isomorphism $\br\, k_v^h\cong\br\, k_v$) is
precisely the one obtained by the snake-lemma construction at the beginning of this
section. It remains to apply Lemma \ref{belem}.
\end{dem}

\section{Proof of Theorem 1}\label{secthm1}

We now prove Theorem \ref{thm1}, of which we take up the notation and assumptions. As in
the case of abelian varieties, the idea is to relate the Brauer-Manin pairing for a torsor
$X$ under a semi-abelian variety $G$ to a Cassels-Tate type pairing. In our case it is the
generalised pairing for 1-motives
$$
\langle\, ,\, \rangle:\, \Sha(M)\times \Sha (M^*)\to\Q/\Z
$$
defined in \cite{hasza}, where $\Sha(M^*)$ is the Tate-Shafarevich group attached to the
dual 1-motive $M^*=[\widehat T\to A^*]$ of $M=[0\to G]$. For the various generalities
about 1-motives used here and in the sequel, we refer to the first section of
\cite{hasza}.

To relate the two pairings, we shall construct a map
$$
\iota:\,\Sha (M^*)\to \Be(X)
$$
and prove that the equality
\begin{equation}\label{beeq}
\langle [X], \beta\rangle =\be (\iota(\beta))
\end{equation}
holds for all $\beta\in \Sha(M^*)$ up to a sign. Theorem \ref{thm1} will then follow from
the non-degeneracy of the Cassels-Tate pairing proven in \cite{hasza}.

To construct the map $\iota$ we proceed as follows. Recall that for a smooth
quasi-projective variety $\overline V$ over a field of characteristic 0 there exists a
generalised Albanese variety $\Alb_V$ introduced in \cite{moruniv} over an algebraically
closed field, and in \cite{rama} in general. It is a semi-abelian variety, and according
to a result of Severi generalised by Serre \cite{moruniv2} and amplified in \cite{rama}
the Cartier dual of the 1-motive $[0\to \Alb_V]$ is $[\div^{\infty,\rm alg}_{V^c}\to
\pic^0_{V^c}]$. Here $V^c$ is a smooth compactification of $V$, and the term
$\div^{\infty,\rm alg}_{V^c}$ is the group of divisors on $V^c$ algebraically equivalent
to 0 and supported in $V^c\setminus V$ viewed as an \'etale locally constant group scheme.
In the case $V=X$ we have $\Alb_V=G$ by definition, and therefore
\begin{equation}\label{m*is}
M^*=[\div^{\infty,\rm alg}_{X^c}\to \pic^0_{X^c}].
\end{equation}
Since there is a natural map of complexes of $k$-group schemes
\begin{equation}\label{mapmdagger}
[\div^{\infty,\rm alg}_{X^c}\to \pic^0_{X^c}]\to [\div^{\infty}_{X^c}\to \pic_{X^c}],
\end{equation}
passing to hypercohomology yields a map
$$
\hyp^1(k, M^*)\to \hyp^1(k, DP(X^c))
$$
with the notation $$DP(X^c):= [\div^\infty_{X^c}\to \pic_{X^c}].$$ Over a number field $k$
the group $\hyp^1(k, DP(X^c))$ is isomorphic to $\br_a\, X$ by Corollary \ref{corbra}, and
the same holds over the completions of $k$. Since the above map is manifestly functorial
for field extensions, we obtain the required map $\iota$ by  restricting to locally
constant elements.

We can thus rewrite the map $\iota$ as
$$
\iota:\, \Sha(M^*)\to\Sha(DP(X^c)).
$$
The previous construction also yields a dual map
\begin{equation}\label{dualmap}
\hyp^1(k, Hom(DP(X^c), \G[1]))\to \hyp^1(k, M)
\end{equation}
by applying the functor $Hom(\quad, \G [1])$ to the map (\ref{mapmdagger}) and taking
hypercohomology (recall that $M=Hom(M^*, \G[1])$). Restricting to locally trivial elements
yields a map
$$
\iota^D:\, \Sha(Hom(DP(X^c), \G[1]))\to \Sha(M).
$$

\begin{rema}\label{remhom}\rm The Hom-functor used in the above formulas is an internal Hom-functor in the bounded derived category of sheaves on the big \'etale site of $\spec k$ restricted to the full subcategory $Sm/k$ of smooth $k$-schemes. The Barsotti-Weil formula $A^*\cong Ext^1(A, \G)$ for abelian schemes holds in this context, because (as O. Wittenberg kindly explained to us) the proof of \cite{oort}, Corollary 17.5 carries over from the {\em fpqc} site to the big \'etale site in the case of smooth group schemes. Hence so does the isomorphism $M^*\cong Hom(M, \G[1])$ used above.  All this remains true if we work instead of $\spec k$ over a suitable open set $U\subset\spec\calo_k$.
\end{rema}

We shall also need versions of the maps constructed above over a suitable
$U\subset\spec\calo_k$. Over $U$ sufficiently small the complex $DP(X^c)$, viewed as a
complex of \'etale sheaves on $Sm/k$, extends to a complex
$$\mathcal{DP}(\mathcal{X}^c):= [\div^\infty_{\mathcal{X}^c/U}\to \pic_{\mathcal{X}^c/U}],$$
where we have used the notations of Remark \ref{braremas} (2). Shrinking $U$ if necessary,
we can also extend the $1$-motive $M$ to a $1$-motive ${\mathcal M}$ over $U$. The main
point is then:

\begin{lem}
The dual 1-motive ${\mathcal M}^*=[{\mathcal Y} \to {\mathcal A}^*]$ is isomorphic to
$[\div^{\infty, {\rm alg}}_{{\mathcal X}^c/U} \to \pic^0_{{\mathcal X}^c/U}]$ over
sufficiently small $U$, where $\div^{\infty, {\rm alg}}_{{\mathcal X}^c/U}$ is the inverse
image of $\pic^0_{{\mathcal X}^c/U}$ in $\div^{\infty}_{{\mathcal X}^c/U}$.
\end{lem}

Here $\pic^0_{{\mathcal X}^c/U}\subset \pic_{{\mathcal X}^c/U}$ is the subsheaf of
elements whose restriction to each fibre of the map ${\mathcal X}^c\to U$ lies in $\pic^0$
of the fibre.\medskip

\begin{dem}
This should be part of a duality theory of Albanese and Picard 1-motives for smooth
quasi-projective schemes over $U$. Since we do not know an adequate reference for this, we
have chosen to circumvent the problem as follows. The group schemes $\pic^0_{{\mathcal
X}^c/U}$ and ${\mathcal A}^*$ are both smooth group schemes of finite type over $U$,
whereas  $\div^{\infty, {\rm alg}}_{{\mathcal X}^c/U}$ and ${\mathcal Y}$ are both
character groups of $U$-tori, so maps defined between their generic fibres extend to maps
over suitable $U$.
\end{dem}

\begin{cor}\label{corid} Over suitable $U\subset\spec\calo_k$ the map $\iota$ lifts to  a map
$$\iota_U:\, \hyp^1_c(U, {\mathcal M}^*)\to \hyp^1_c(U, {\mathcal DP}({\mathcal X}^c)),$$ and the map
 (\ref{dualmap})  extends to
 $$\iota^D_U:\, \hyp^1(U, Hom({\mathcal DP}(\mathcal{X}^c), \G[1]))\to \hyp^1(U, {\mathcal M}).$$
\end{cor}

\begin{dem} The map (\ref{mapmdagger}) extends to a map of complexes
$$[\div^{\infty, {\rm alg}}_{{\mathcal X}^c/U} \to \pic^0_{{\mathcal X}^c/U}] \to
{\mathcal DP}({\mathcal X}^c),$$ so by the lemma we dispose of a map ${\mathcal
M}^*\to{\mathcal DP}({\mathcal X}^c)$. The required maps are obtained by passing to
cohomology.
\end{dem}

In the previous section we worked with a certain extension class ${\mathcal E}_X$ in
$\Ext^1_{Sm/U}(\mathcal{KD'(X)}, \G[1])$. According to Remark \ref{braremas} (2) the
Ext-group here is isomorphic to $\Ext^1_{Sm/U}(\mathcal{DP}({\mathcal X}^c), \G[1])$. The
next lemma will imply that we may identify ${\mathcal E}_X$ with a class $\mathcal{E_X'}$
in the group $\hyp^1(U, Hom(\mathcal{DP}({\mathcal X}^c), \G[1])$, i.e. the source of the
map $\iota^D_U$.

\begin{lem}\label{extiso} There are canonical isomorphisms
$$\Ext^j_{Sm/U}(\mathcal{DP}({\mathcal X}^c), \G[1])\cong \hyp^j(U,
Hom(\mathcal{DP}({\mathcal X}^c), \G[1]))$$ for all $j>0$. \end{lem}

\begin{dem} We start with the isomorphism
$$
\Ext^j_{Sm/U}(\mathcal{DP}({\mathcal X}^c), \G[1])\cong \hyp^j(U, {\bf
R}Hom_{Sm/U}(\mathcal{DP}({\mathcal X}^c), \G[1]))
$$
coming from the derived category analogue of the spectral sequence for composite functors;
here ${\bf R}Hom$ denotes the derived functor of the inner Hom explained in Remark
\ref{remhom}. It shows that the lemma follows from the vanishing of the groups
$$Ext_{Sm/U}^{i-1}(\mathcal{DP}({\mathcal X}^c),
\G[1])=Ext^i_{Sm/U}(\mathcal{DP}({\mathcal X}^c), \G)$$ for $i>1$ and $i=0$. This may be
checked on geometric fibres, so we are reduced to checking
$\Ext^i_{Sm/\kbar}(DP(X^c)_\kbar, \G)=0$ for $i>1$ and $i=0$, over an arbitrary
algebraically closed field $\kbar$. We drop the subscripts in the rest of the proof.

Observe first that the cokernel of the map of complexes (\ref{mapmdagger}) is
quasi-isomorphic to the complex $[0\to B(X^c)]$, where $B(X^c)$ is the quotient of the
N\'eron-Severi-group of $X^c$ modulo the subgroup of classes coming from divisors at
infinity; in particular, its $\kbar$-points form a finitely generated abelian group. Hence
the sheaf  $\Ext^i(B(X^c), \G)$ is trivial for $i>0$ (see \cite{skobook}, Sublemma 2.3.8).
Therefore the distinguished triangle coming from   (\ref{mapmdagger}) shows that it is
enough to prove $\Ext^i([\div^{\infty,\rm alg}_{X^c}\to \pic^0_{X^c}],\G)=0$ for $i>1$ and
$i=0$, which is the same as proving $\Ext^i(M^*, \G)=0$ by the isomorphism (\ref{m*is}).
The case $i=0$ then follows from the fact that every morphism from $A^*$ to $\G$ is
trivial. For the case $i>1$ we remark that the stupid filtration on $M^*=[\widehat T\to
A^*]$ induces an exact sequence
$$
\Ext^{i}(A^*, \G)\to {\Ext}^i(M^*, \G)\to {\Ext}^{i-1}(\widehat T, \G).
$$
Here the terms at the two extremities are trivial for $i>1$ (the left one by \cite{oort},
Prop. 12.3), hence so is the middle one.
\end{dem}

\begin{remas}\label{corextiso}\rm ${}$\smallskip

\noindent 1. Here it was crucial to work with extensions over the big \'etale site; over
the small \'etale site of $\kbar$ the group $\Ext^{i}(A^*, \G)$ is trivial even for $i=1$.
\smallskip

\noindent 2. In the course of the proof we also established canonical isomorphisms
$$\Ext^j_{Sm/U}({\mathcal M}^*, \G[1])\cong \hyp^j(U,
Hom({\mathcal M}^*, \G[1]))$$ for all $j>0$, that in turn yield isomorphisms
$$\Ext^j_{Sm/k}({M}^*, \G[1])\cong \hyp^j(k, Hom({ M}^*, \G[1]))$$ at the generic point.
\end{remas}

By the lemma we may apply the map $\iota^D_U$ to ${\mathcal E}_X'$, and obtain a class in
$\hyp^1(U, {\mathcal M})$. Over the generic point the latter group becomes $H^1(k, G)$,
and denoting by $E_X'$ the image of ${\mathcal E}_X'$ in $\hyp^1(k, Hom(DP(X^c), \G[1])$
we have the following basic fact.

\begin{lem}\label{44} The image of $E_X'$ by the map (\ref{dualmap}) equals (up to a sign) the class of $X$ in $H^1(k,G)=\hyp^1(k,M)$.
\end{lem}

The lemma should be true over an arbitrary field of characteristic 0. It is known in the
two extreme cases $G=A$ and $G=T$ (see references in the proof below); we leave the
general case to the reader as a challenge. The following proof, which is sufficient for our
purposes, works under the assumptions of Theorem \ref{thm1} (i.e. over a number field,
assuming $X(\A_k)\neq\emptyset$ and the finiteness of $\Sha(A)$). Also, as O. Wittenberg
pointed out to us, Corollary 4.2.4 of \cite{witten} implies that $E_X'$ maps to 0 in
$H^1(k, G)$ if and only if $[X]=0$, which is also sufficient for the proof of Theorem
\ref{thm1} given below.\medskip

\noindent{\em Proof of Lemma \ref{44}}. Thanks to Proposition 2.1 of \cite{russian} the
case $G=A$ is known, and we may complete the proof of Theorem \ref{thm1} given below in
this special case. Thus we are allowed to apply Theorem \ref{thm1} to the pushforward
torsor $p_*X$ under $A$ (here of course $p:\, G\to A$ is the natural projection map), and
conclude that it is trivial. The exact sequence
$$
H^1(k, T)\stackrel{i_*}\to H^ 1(k, G)\stackrel{p_*}\to H^1(k, A)
$$
then implies that $X=i_*Y$ for some $k$-torsor $Y$ under $T$, where ${i: T\to G}$ is the
natural inclusion.

The map (\ref{dualmap}) factors through  $\hyp^1(k, Hom(M^*, \G[1]))$ by construction, and
by Remark \ref{corextiso} (2) we may identify the image of $E_X'$ in the latter group with
a class $E_X^0\in \Ext^1(M^*, \G[1])$. By performing the same construction for the torsor
$Y$ we obtain a class $E_Y^0\in \Ext^1(\widehat T[1], \G[1])$. According to \cite{witten},
Proposition 4.1.4 applied with $V=Y$ and $W=X$ we have $E_X^0=i^*E_Y^0$, where $i^*:\,
\Ext^1(\widehat T[1], \G[1])\to \Ext^1(M^*, \G[1])$ is the natural map induced by the
projection $M^*\to \widehat T[1]$. (Note that this equality is not completely obvious,
because the map $Y\to X$ is not dominating.) But for $G=T$ the lemma is known over an
arbitrary field (\cite{skobook}, Lemma 2.4.3), so the image of $E_Y^0$ in $H^1(k, T)$  is
$[Y]$ up to a sign. The image of $[Y]$ in $H^1(k, G)$ is $[X]$, so the lemma in the
general case follows from the commutativity of the diagram
$$
\begin{CD}
H^1(k, T) @>>> H^1(k, G) \\
@AAA @AAA \\
\Ext^1(\widehat T[1], \G[1]) @>>> \Ext^1(M^*, \G[1]).
\end{CD}
$$

\enddem

\noindent{\em Proof of Theorem \ref{thm1}}. As already remarked at the beginning of this
section, for the proof of the theorem it will suffice to verify formula (\ref{beeq}) for
the class of $X$ in $\Sha(M)$, i.e. the equality $\langle [X], \beta\rangle =\be
(\iota(\beta))$ up to a sign for all $\beta\in\Sha(M^*)$. Indeed, our assumption that $X$
has an adelic point orthogonal to $\Be(X)$ implies the triviality of the map $\be$, so the
right hand side of the formula is 0 for all $\beta$ in $\Sha(M^*)$.  But under the
finiteness assumption on $\Sha(A)$ the Cassels-Tate pairing $\langle\,,\,\rangle$ is
non-degenerate (\cite{hasza}, Corollary 4.9), so $[X]=0$, i.e. $X$ has a $k$-rational
point.

We now verify formula (\ref{beeq}). Consider the cup-product pairing
$$
\hyp^1(U, Hom(\mathcal{DP}(X^c), \G[1]))\times \hyp^1_c(U, \mathcal{DP}(X^c))\to H^3_c(U,
\G)
$$
and recall that we have defined above a class $\mathcal{E}_X'$ in the cohomology group
$\hyp^1(U,Hom(\mathcal{DP}(X^c), \G[1]))$. By construction, taking the product of the
class $\mathcal{E}_X'$ with some $\alpha_U\in\hyp^1_c(U, \mathcal{DP}(X^c))$ under this
pairing is the same as the element ${\mathcal E}_X\cup \alpha_U$ considered in Proposition
\ref{BM}. So applying the proposition we obtain the equality ${\mathcal E}_X'\cup
\alpha_U=\be(\alpha)$ in the case when $\alpha_U$ maps to a locally trivial element in
$\hyp^1(k, DP(X^c))$.

Moreover, using the maps constructed in Corollary \ref{corid} we have a diagram

\vbox{$$ \hyp^1(U, \mathcal{M})\qquad\qquad  \times\qquad \qquad\hyp^1_c(U, \mathcal{M}^*)
\quad\to\quad H^3_c(U,\G)
$$

$$
\iota^D_U\uparrow\qquad\qquad\qquad\qquad \qquad\qquad
\downarrow\quad\iota_U\quad\qquad\qquad \qquad\downarrow\, {\rm id}
$$

$$
\hyp^1(U, Hom(\mathcal{DP}(X^c), \G[1]))\times \hyp^1_c(U, \mathcal{DP}(X^c))\to H^3_c(U,
\G)
$$}

\noindent where the horizontal maps are cup-product pairings. The diagram commutes by
construction, and the image $\iota^D(E_X')$ of the element $\iota^D_U(\mathcal{E}_X')$ in
$H^1(k, M)=H^1(k, G)$ is the class $[X]$ up to a sign by the previous lemma. By Corollary
4.3 of \cite{hasza} each element $\beta\in\Sha(M^*)$ comes from some $\beta'\in\hyp^1_c(U,
\mathcal{M}^*)$ for $U$ sufficiently small, and moreover the value of the upper pairing on
$(\iota^D_U(\mathcal{E}_X'), \beta')$ equals the value of the Cassels-Tate pairing on
$(\iota^D({E}_X'), \beta)$, i.e. on $([X], \beta)$ up to a sign. The commutativity of the
diagram together with the arguments of the previous paragraph implies that this value
equals $\be(\iota(\beta))$. This proves formula (\ref{beeq}), and thereby the theorem.
\enddem

\begin{rema}\label{berema}\rm As a complement to the theorem, we justify here a claim made in the introduction, namely that the group $\Be (X)$ is {\em finite}. In \cite{bcs}, Proposition 2.14 the finiteness of $\Be (V)$ is verified for a smooth {\em proper} $V$ such that $\Sha(\pic^0\, V)$ is finite. To deduce the statement for our $X$, apply this result with $V=X^c$, a smooth compactification of $X$. The condition on the Tate-Shafarevich group holds because $\pic^0(V)$ is the Picard variety of the Albanese variety of $V$ (theorem of Severi), and the latter is none but $A$ (see \cite{moruniv}, \cite{moruniv2} for these facts). It remains to add that $\Be(V)\cong\Be(X)$ in view of \cite{sansuc}, (6.1.4).
\end{rema}

To conclude this section we mention a variant of Theorem~\ref{thm1} that deals with points
of $X$ over the direct product $k_\Omega$ of all completions of $k$ instead of $\A_k$. In
this situation we look at a modified version of the Brauer-Manin pairing, namely the
induced pairing
\begin{equation}\label{nrpairing}
X(k_\Omega)\times (\brnr X/\br\, k) \to {\bf Q}/{\bf Z},
\end{equation}
where $\brnr X$ is the unramified Brauer group of $X$, which may be defined as the Brauer
group of a smooth compactification $V$ of $X$. Since $\Be(X)\cong \Be(V)$ as in the remark
above, the group $\Be(X)$ is contained in $\brnr X/\br\, k$.

\begin{cor}
Let $G$ be a semi-abelian variety defined over $k$, and let $X$ be a $k$-torsor under $G$.
Assume that the Tate-Shafarevich group of the abelian quotient of $G$ is finite. If there
is a point of $X(k_{\Omega})$ annihilated by all elements of $\Be (X)$ under the pairing
(\ref{nrpairing}), then $X$ has a $k$-rational point.
\end{cor}

The corollary immediately follows from Theorem \ref{thm1} and the following lemma:

\begin{lem}
Let $X$ be a smooth geometrically integral variety defined over $k$.
If there is a point of $X(k_{\Omega})$ orthogonal to $\Be(X)$ under the pairing
(\ref{nrpairing}), then there is also an adelic point on $X$ orthogonal to $\Be(X)$ under
the pairing (\ref{bepairing}).
\end{lem}

\begin{dem} From Chow's lemma we know that $X$ contains a quasi-projective
  open subset $U$. 
Choose a finite set $S$ of places of $k$ such that the pair $U\subset X$
extends to a pair of smooth
schemes ${\mathcal U}\subset{\mathcal X}$ over $\spec (\calo_{k,S})$ with $U$ quasi-projective, where $\calo_{k,S}$ is
the ring of $S$-integers of $k$. From the Lang--Weil estimates and Hensel's
lemma we know that by enlarging $S$ if necessary we have ${\mathcal U}(\calo_v) \neq
\emptyset$ for $v \not \in S $, and hence the same holds for $\mathcal X$. Now if
$(P_v) \in X(k_{\Omega})$ is orthogonal to $\Be(X)$, we replace $P_v$ by an
$\calo_v$-point $P'_v$ of ${\mathcal X}$ for $v \notin S$. Then $(P'_v)$ is an adelic
point of $X$, and this adelic point remains orthogonal to $\Be(X)$ because elements of
$\Be(X)$ induce constant elements of $\br (X \times_k k_v)$ for every place $v$.
\end{dem}

\section{The Cassels-Tate dual exact sequence for 1-motives}\label{catate}

In this section we prove Theorem \ref{casstatetheo}, of which we take up the notation.
Recall that by convention for an archimedean place $v$ of $k$ the notation $\hyp^0(k_v,M)$
stands for the Tate group $\widehat \hyp^0(k_v,M)$, which is a $2$-torsion finite group.
Also, recall from (\cite{hasza}, \S 2) that the group $\hyp^0(k_v,M)$ is equipped with a
natural topology. In the case $M=G$ and $v$ finite, this is just the usual $v$-adic
topology on $H^0(k_v,G)=G(k_v)$, but in general the topology on $\hyp^0(k_v,M)$ is not
Hausdorff.

We denote by $\overline{{\hyp^0(k,M)}}$ the closure of the diagonal image of $\hyp^0(k,M)$
in the topological direct product of the $\hyp^0(k_v,M)$. The local pairings
$(\,\,,\,\,)_v$ of (\cite{hasza}, \S 2) induce a map $$\theta : \prod_{v\in\Omega}
\hyp^0(k_v,M) \to \Sha^1_{\omega}(M^*)^D$$ defined by
$$\theta((m_v))(\alpha)=\sum_{v \in \Omega} (m_v,\alpha_v)_v,$$
where $\alpha_v$ is the image of $\alpha \in \Sha^1_{\omega}(M^*)$ in $\hyp^1(k_v,M)$ (the
sum is finite by definition of $\Sha^1_{\omega}(M^*)$). On the other hand, the analogue of
Cassels-Tate pairing for 1-motives (\cite{hasza}, Theorem 4.8) and the inclusion
$\Sha^1(M^*) \subset \Sha^1_{\omega}(M^*)$ induce a map
$$p : \Sha^1_{\omega}(M^*)^D \to \Sha^1(M)$$
We have thus defined all maps in the sequence
$$0 \to \overline{{\hyp^0(k,M)}} \to \prod_{v \in \Omega} \hyp^0(k_v,M)
\stackrel{\theta}{\to} \Sha^1_{\omega}(M^*)^D \stackrel{p}{\to} \Sha^1(M) \to 0$$ and our
task is to prove its exactness.

We shall need several intermediate results. The first one is the following well-known
lemma, for which we give a proof by lack of a reference.

\begin{lem}
Let $Y$ be a $k$-group scheme \'etale locally isomorphic to $\Z^r$ for some $r>0$. Then
the group $\Sha^2_{\omega} (Y)$ is finite.
\end{lem}

Here by definition $\Sha^2_\omega(Y):=\Sha^1_\omega([Y\to 0])$, with the notation of the
introduction.

\begin{proof} Let $L$ be a finite Galois extension of $k$ that splits $Y$. Since $\Sha^2_{\omega}(\Z)=\Sha^1_{\omega}(\Q/ \Z)$ is
zero by Chebotarev's density theorem, we obtain that $\Sha^2_{\omega} (Y)$ is a subgroup
of $H^2(\gal(L|k),Y)$, which is a torsion group annihilated by $n=[L:k]$. The boundary map
$$H^1(\gal(L|k),Y/n Y) \to H^2(\gal(L|k),Y)$$ obtained from the exact sequence of
$\gal(L|k)$-modules
$$0 \to Y \to Y \to Y/n Y \to 0$$ is therefore surjective. Since $\gal(L|k)$ and $Y/nY$ are finite, the lemma follows.
\end{proof}

Now return to the situation above, and recall from \cite{hasza}, Theorem 2.3 and Remark
2.4 that the local pairings $(\,\,,\,\,)_v$ used in the definition of $\theta$ actually
factor through the profinite completion $\hyp^0(k_v,M)^{\wedge}$ of
$\hyp^0(k_v,M)$, hence $\theta$ extends to $\hyp^0(k_v,M)^{\wedge}$.
Technical
complications will arise from the fact that the topology on $\hyp^0(k_v,M)$ is in general
finer than the topology induced by the profinite topology of $\hyp^0(k_v,M)^{\wedge}$. For
instance, this is the case for $M=[0\to T]$ with $T$ a torus.

\begin{lem} \label{nocompl}
The groups $\prod_{v \in \Omega} \hyp^0(k_v,M)^{\wedge}$ and $\prod_{v \in \Omega}
\hyp^0(k_v,M)$ have the same image by $\theta$.
\end{lem}

\begin{dem} For $v$ archimedean, the group
$\hyp^0(k_v,M)$ is finite, hence it is the same as its profinite completion, so we can
concentrate on the finite places. We proceed by d\'evissage, starting with the case
$M=[0\to G]$. Let $v$ be a finite place of $k$. Since $A$ is proper, we have
$H^0(k_v,A)=A(k_v)=H^0(k_v,A)^{\wedge}$. Using the exact sequences
$$0 \to T(k_v) \to G(k_v) \to A(k_v) \to H^1(k_v,T)$$
$$0 \to T(k_v)^{\wedge} \to G(k_v)^{\wedge} \to A(k_v) \to H^1(k_v,T)$$
(cf. \cite{hasza}, Lemma
2.2), we obtain that
$$\prod_{v \in \Omega} H^0(k_v,G)^{\wedge}=\left\{g+t:\, g \in \im\left(\prod_{v \in \Omega} H^0(k_v,G)\right),\, t \in
\prod_{v \in \Omega} H^0(k_v,T)^{\wedge}\right\}.$$ Therefore it is sufficient to prove
the statement for $G=T$. But this follows from the facts that $\Sha^1_\omega
(M^*)=\Sha^2_{\omega}(Y^*)$ is finite (by the previous lemma), and each $H^0(k_v,T)$ is
dense in $H^0(k_v,T)^{\wedge}$.

\smallskip

The same method reduces the general case to the case $M=[0\to G]$, using the exact
sequences (\cite{hasza}, p. 101):
$$H^0(k_v,G) \to \hyp^0(k_v,M) \to H^1(k_v,Y) \to H^1(k_v,G)$$
$$H^0(k_v,G)^{\wedge} \to \hyp^0(k_v,M)^{\wedge} \to H^1(k_v,Y) \to H^1(k_v,G)$$
\end{dem}

Denote by $\Sha^1_S(M^*)$ the kernel of the diagonal map $$\hyp^1(k,M^*) \to \prod_{v \not
\in S} \hyp^1(k_v,M^*).$$ As above, the local pairings induce maps 
$$\theta_S : \prod_{v\in S} \hyp^0(k_v,M) \to \Sha^1_S(M^*)^D$$ and
$$\widehat\theta_S : \prod_{v\in S} \hyp^0(k_v,M)^\wedge \to \Sha^1_S(M^*)^D.$$

\begin{prop}
Let $S$ be a finite set of places of $k$.  Assume that the Tate-Shafarevich group $\Sha
  (A)$ of the abelian quotient of $G$ is finite.
\begin{enumerate}
\item The sequence
\begin{equation}\label{kacsa}
\hyp^0(k,M)^{\wedge} \to \prod_{v \in S} \hyp^0(k_v,M)^{\wedge}
\stackrel{\widehat\theta_S}{\to} \Sha^1_S(M^*)^D
\end{equation}
 is exact.
\item Denote by $\overline{{\hyp^0(k,M)}}_S$ the closure of the diagonal image of
$\hyp^0(k,M)$ in $\prod_{v \in S} \hyp^0(k_v,M)$. Then the sequence
\begin{equation}\label{finisequ}
0 \to \overline{{\hyp^0(k,M)}}_S \to \prod_{v \in S} \hyp^0(k_v,M)
\stackrel{\theta_S}{\to} \Sha^1_S(M^*)^D
\end{equation}
is exact.
\end{enumerate}
\end{prop}

\begin{dem} 1. Let $\P^1(M^*)$ be the restricted product of the $\hyp^1(k_v,M^*)$
(cf. \cite{hasza}, \S 5). By the Poitou-Tate exact sequence for 1-motives (\cite{hasza},
Th. 5.6), there is an exact sequence
$$\hyp^1(k,M^*) \to \P^1(M^*)_{\rm tors} \to (\hyp^0(k,M)^D)_{\rm tors}.$$
(Recall that this uses the finiteness of the Tate-Shafarevich group of $A^*$, which is
equivalent to that of $A$ by \cite{adt}, Remark I.6.14(c)). Sending an element of
$\prod_{v \in S} \hyp^1(k_v,M^*)$ to $\P^1(M^*)_{\rm tors}$ via the map
$$(m_v)_{v \in S} \mapsto ((m_v),0,0,...)$$ yields an exact sequence of discrete torsion
groups
$$\Sha^1_S(M^*) \to \prod_{v \in S} \hyp^1(k_v,M^*) \to
(\hyp^0(k,M)^D)_{\rm tors}.$$ We claim that the required exact sequence is the dual of the
above. Indeed, the dual of the discrete torsion group $\hyp^1(k_v,M^*)$ is the profinite
group $\hyp^0(k_v,M)^{\wedge}$ by the local duality theorem (\cite{hasza}, Th. 2.3), and
the dual of the discrete torsion group $(\hyp^0(k,M)^D)_{\rm tors}$ is the profinite
completion $\hyp^0(k,M)^{\wedge}$ of $\hyp^0(k,M)$, because $(\hyp^0(k,M)^D)_{\rm tors}$
is nothing but the direct limit (over the subgroups $I \subset \hyp^0(k,M)$ of finite
index) of the groups $\Hom(\hyp^0(k,M)/I,\Q/\Z)$.

\smallskip

2. Consider the commutative diagram
\begin{equation} \label{second}
\begin{CD}
\hyp^0(k,M) @>j>> \prod_{v \in S} \hyp^0(k_v,M)@>{\theta_S}>> \Sha^1_S(M^*)^D \cr @VVV
@VVV @VV{\rm id}V \cr \hyp^0(k,M)^{\wedge} @>{j^\wedge}>> \prod_{v \in S}
\hyp^0(k_v,M)^{\wedge} @>{\widehat\theta_S}>> \Sha^1_S(M^*)^D
\end{CD}
\end{equation}
The second line is exact by what we have just proven, so the first line is a complex.
Hence so is the sequence (\ref{finisequ}) by continuity of $\widehat\theta_S$. Denote by $J$ the
closure of the image of $j$ in the above diagram.  Set $$C:=\prod_{v \in S}
\hyp^0(k_v,M)/J,$$ and equip $C$ with the quotient topology. In particular, $C$ is a
Hausdorff topological group (because $J$ is closed).

Consider now the commutative diagram
\begin{equation}\label{third}
\begin{CD}
0 @>>> J @>>> \prod_{v \in S} \hyp^0(k_v,M) @>>> C \cr && @VVV @VVV @VVV \cr && J^{\wedge}
@>>> \prod_{v \in S} \hyp^0(k_v,M)^{\wedge} @>>> C^{\wedge}.
\end{CD}
\end{equation}

Assume for the moment that the right vertical map here is injective. We can
then derive
the exactness of sequence (\ref{kacsa}) as follows. The first line of diagram (\ref{third}) is exact by definition, and the second line is a complex
because it is the completion of an exact sequence. 
Since the second line of diagram (\ref{second}) is exact, the image of an
element ${x\in\ker(\theta_S)}$ in $\ker(\widehat\theta_S)$
comes from $H^0(k,M)^{\wedge}$, 
hence from $J^{\wedge}$. A diagram chase in (\ref{third}) then shows
$x\in J$, which is what we wanted to prove.

Now the injectivity of the right vertical map in (\ref{third}) follows from statement (3) of the Appendix to \cite{hasza}, of which we have
to check the assumptions. The last horizontal map above is an open mapping because it is a
quotient map. The group $C$ is Hausdorff, locally compact and totally disconnected by
construction; it remains to check that it is also compactly generated (i.e. it is
generated as a group by the elements of a compact subset). This is because by \S 2 of
\cite{hasza} the group $\hyp^0(k_v,M)$ has a finite index open subgroup that is a
topological quotient of $H^0(k_v,G)$, so $C$ has a finite index open subgroup $C'$ that is
a quotient of the product of the $H^0(k_v,G)$ for $v\in S$. Since each $H^0(k_v,G)$ is
compactly generated (this follows from the theory of $p$-adic Lie groups) and $C'$ is
Hausdorff, we obtain that $C'$, and hence $C$, are compactly generated.

\end{dem}

\noindent{\it Proof of Theorem~\ref{casstatetheo}.} Let us start by proving the exactness
of the sequence
\begin{equation} \label{firstsequ}
\prod_{v \in \Omega} \hyp^0(k_v,M) \to \Sha^1_{\omega}(M^*)^D \to \Sha^1(M) \to 0.
\end{equation}

The sequence
\begin{equation} \label{dualsequ}
0 \to \Sha^1(M^*) \to \Sha^1_{\omega}(M^*) \to \bigoplus_{v \in \Omega} \hyp^1(k_v,M^*)
\end{equation}
is exact by definition. By the local duality theorem for 1-motives (\cite{hasza}, Theorem
2.3 and Proposition 2.9), the dual of each group $\hyp^1(k_v,M^*)$ is
$\hyp^0(k_v,M)^{\wedge}$, and by the global duality theorem (\cite{hasza}, Corollary 4.9),
the dual of $\Sha^1(M^*)$ is $\Sha^1(M)$ under our finiteness assumption on
Tate-Shafarevich groups. Therefore the dual of (\ref{dualsequ}) is the exact sequence
$$\prod_{v \in \Omega} \hyp^0(k_v,M)^{\wedge} \to \Sha^1_{\omega}(M^*)^D
\to \Sha^1(M) \to 0,$$ and Lemma~\ref{nocompl} gives the exactness of (\ref{firstsequ}).

\smallskip

It remains to prove the exactness of the sequence
$$0 \to \overline{{\hyp^0(k,M)}} \to \prod_{v \in \Omega}
\hyp^0(k_v,M) \to \Sha^1_{\omega}(M^*)^D.$$ But this sequence is obtained by applying the
(left exact) inverse limit functor (over all finite subsets $S \subset \Omega$) to the
exact sequences (\ref{finisequ}). Indeed, by definition of the direct product topology the
inverse limit of the groups $\overline{{\hyp^0(k,M)}}_S$ is $\overline{{\hyp^0(k,M)}}$,
and the inverse limit of the groups $\Sha^1_S(M^*)^D$ is the dual of the direct limit of
the discrete torsion groups $\Sha^1_S(M^*)$, i.e. the dual of $\Sha^1_{\omega}(M^*)$.
\enddem

\section{Obstruction to weak approximation}\label{secaf}

In the study of weak approximation on a variety $X$ one works with a modified version of
the Brauer-Manin pairing, namely with the induced pairing
$$
X(k_\Omega)\times \brnr X \to {\bf Q}/{\bf Z},
$$
already encountered at the end of Section \ref{secthm1}, where $k_\Omega$ is the
topological direct product of all completions of $k$, and $\brnr X$ is the unramified
Brauer group of $X$. One may also work with subgroups of $\brnr X$, such as $\brnr_1\,
X:=\ker(\brnr X \to\brnr (X\times_k\kbar))$. Finally, for a smooth $k$-group scheme $G$
there is yet another variant, which is the one we shall use:
\begin{equation}\label{pairing2}
\prod_{v\in\Omega} H^0(k_v,G) \times \brnr_1\,G\to {\bf Q}/{\bf Z}.
\end{equation}
Here we have taken the same convention at the archimedean places as in  Theorem
\ref{casstatetheo} proven above. Concerning this pairing one has the following result,
first proven in \cite{imrn}:

\begin{theo}\label{thm2} Let $G$ be a semi-abelian variety defined over $k$. Assuming that the abelian quotient has finite Tate-Shafarevich group, the left kernel of the pairing
(\ref{pairing2}) is contained in the closure of the diagonal image of ${{G(k)}}$.
\end{theo}

This result was proven in {\em loc. cit.} for arbitrary connected algebraic groups, but
the key case is that of a semi-abelian variety. We now show that the statement can be
easily derived from Theorem \ref{casstatetheo} as follows. The Brauer-Manin pairing
induces a map
$$\prod_{v\in\Omega} H^0(k_v, G) \to (\brnr_1 G/\br \, k)^D.$$ Going through the construction of the map $\iota$  at the beginning of Section \ref{secthm1} with $\Sha_\omega$ in place of $\Sha$ we obtain a map
$\iota_\omega:\,\Sha^1_{\omega}(M^*) \to \Be_\omega (G)$, where $\Be_\omega(G)\subset
\br_a\, G$ is the subgroup of elements that are locally trivial for almost all places for
$k$. Using the inclusion $\Be_\omega(X)\subset\brnr X/\br\, k$ resulting from
(\cite{sansuc}, 6.1.4) we thus obtain a map $r : \Sha^1_{\omega}(M^*) \to \brnr_1\, G
/\br\, k$, whence a diagram
$$
\begin{diagram}
\node{}\node{} \node{(\brnr_1 G /\br\, k)^D}\arrow{s,r}{\iota_\omega^D} \\
\node{\overline{{G(k)}}} \arrow{e} \node{\prod_{v\in\Omega} H^0(k_v, G_v)}\arrow{e}
\arrow{ne} \node{\Sha^1_{\omega}(M^*)^D}
\end{diagram}
$$
where $G_v:=G\times_kk_v$. If we prove that the triangle commutes, the theorem follows,
since the bottom row is exact by Theorem~\ref{casstatetheo}.

We shall prove the commutativity of the dualized diagram
\begin{equation}\label{triangle}
\begin{diagram}
\node{\Sha^1_{\omega}(M^*)}\arrow{s,l}{\iota_\omega}\arrow{e}
\node{H^0(k_v, G_v)^D}\\
\node{\brnr_1 G /\br\, k}\arrow{ne}
\end{diagram}
\end{equation}
for all places $v$. Here the horizontal map is induced by local duality, so it is in fact
enough to consider the finitely many nonzero images of an element in
$\Sha^1_{\omega}(M^*)$ by the restriction maps $\Sha^1_{\omega}(M^*)\to H^1(k_v, M^*)$ and
show that the diagram
$$
\begin{diagram}
\node{\hyp^1(k_v, M^*)}\arrow{s}\arrow{e} \node{H^0(k_v, G_v)^D}  \\
\node{\br_a\, G_v} \arrow{ne}
\end{diagram}
$$
commutes, where the diagonal map is induced by the evaluation pairing
\begin{equation}\label{loceval}
G(k_v)\times \br\, G_v\to \br\, k_v\cong\Q/\Z,
\end{equation}
and we view $\br_a\, G_v$ as a subgroup of $\br\, G_v$ thanks to the splitting of the map
$\br\, k_v\to\br\, G_v$ coming from the zero section of $G_v$.

To do so, return to the beginning of Section \ref{secbm} and observe that in the case
$X=G$ the maps (\ref{split}) actually assemble to a pairing of complexes of Galois modules
$$
[0\to G(\kbar_v)]\,\otimes_{\Z}\,[\kbar_v (G)^\times\to \div (G\times_k\kbar_v)]\to
[\kbar_v^\times\to 0].
$$
The sections $\kbar_v (G)^\times\to \kbar_v^\times$ used in this construction are not
canonical, but the pairing becomes canonical at the level of the derived category, again
by the argument in (\cite{skobook}, Theorem 2.3.4 (b)). We thus obtain a cup-product
pairing
$$
H^0(k_v, G)\times \hyp^1(k_v, [\kbar_v (G)^\times\to \div (G\times_k\kbar_v)])\to \br\,
k_v
$$
that identifies via Lemma \ref{bra} with the restriction of the local pairing
(\ref{loceval}) to $\br_1\, G_v$ by the argument at the beginning of Section \ref{secbm}.
On the other hand, we may lift the map $\hyp^1(k_v, M^*)\to \br_a\, G_v$ to a map
$\hyp^1(k_v, M^*)\to \br_1\, G_v$ via the zero section as above, and then the claim
follows from the commutativity of the diagram of cup-product pairings

\vbox{$$ \hyp^0(k_v, M)\times\hyp^1(k_v, M^*)\to \br\, k_v
$$

$$
{\rm id}\quad\downarrow\qquad\qquad\quad\downarrow\quad\qquad\qquad\downarrow\, {\rm id}
$$

$$
H^0(k_v, G)\quad\times\quad \br_1\, G_v\,\,\to \br\, k_v.
$$}

\end{document}